\documentclass{amsart}

\usepackage{graphicx}
\usepackage{mathrsfs, mathscinet, amssymb}
\usepackage{mathtools}
\usepackage{hyperref}

\newtheorem{thm}{Theorem}[section]
\newtheorem{lem}[thm]{Lemma}
\newtheorem{cor}[thm]{Corollary}

\theoremstyle{definition}
\newtheorem{exm}[thm]{Example}
\newtheorem{defn}[thm]{Definition}

\newcommand{\Z}{\mathbb{Z}}
\newcommand{\N}{\Z_{\ge 0}}   % Nonnegative integers

\newcommand{\R}{\mathbb{R}}
\newcommand{\C}{\mathbb{C}}
\newcommand{\gl}{\mathfrak{gl}}
\renewcommand{\emptyset}{\varnothing}
\renewcommand{\epsilon}{\varepsilon}
\newcommand{\gtp}[1]{\mathbf{#1}}
\newcommand{\domby}{\trianglelefteq}
\newcommand{\K}{\mathcal{K}}    % Stretched Kostka coefficient
\newcommand{\primby}{\vartriangleleft}
\newcommand{\PP}{\mathscr{T}}   % Tiling
\newcommand{\A}{A}              % Tiling matrix
\renewcommand{\bar}[1]{\smash{\overline{#1}}}
\newcommand{\prt}{\bar}         % Reordered as partition

\DeclarePairedDelimiter{\size}{\lvert}{\rvert}
\DeclarePairedDelimiter{\card}{\lvert}{\rvert}
\DeclarePairedDelimiter{\setof}{\lbrace}{\rbrace}
\DeclarePairedDelimiter{\opint}{\rbrack}{\lbrack}

\DeclareMathOperator{\hwt}{\mathbf{hwt}}
\DeclareMathOperator{\wt}{\mathbf{wt}}
\DeclareMathOperator{\intr}{int}
\DeclareMathOperator{\Dom}{D}   % Set dominated by argument

\title{%
   Degrees of stretched Kostka coefficients
}%
  
\author{%
   Tyrrell B. McAllister
}%
\address{%
   Department of Mathematics and Computer Science \\
   Eindhoven University of Technology \\
   P.O.~Box 513 \\
   5600 MB Eindhoven \\
   The Netherlands%
}%
\email{%
   tmcallis@win.tue.nl
}%
\thanks{%
   Research supported by NSF VIGRE Grant No.~DMS-0135345 and by
   NWO Mathematics Cluster DIAMANT%
}%



\begin{document}
   \begin{abstract}
      Given a partition \( \lambda \) and a composition \( \beta
      \), the \emph{stretched Kostka coefficient} \( \K_{\lambda
      \beta}(n) \) is the map \( n \mapsto K_{n\lambda, n\beta} \)
      sending each positive integer \( n \) to the Kostka
      coefficient indexed by \( n\lambda \) and \( n\beta \).
      Kirillov and Reshetikhin~\cite{KR86} have shown that
      stretched Kostka coefficients are polynomial functions of \(
      n \).  King, Tollu, and Toumazet have conjectured that these
      polynomials always have nonnegative
      coefficients~\cite{KTT04}, and they have given a conjectural
      expression for their degrees~\cite{KTT06}.
      
      We prove the values conjectured by King, Tollu, and Toumazet
      for the degrees of stretched Kostka coefficients.  Our proof
      depends upon the polyhedral geometry of Gelfand--Tsetlin
      polytopes and uses tilings of GT-patterns, a combinatorial
      structure introduced in~\cite{DLM04}.
   \end{abstract}

   \maketitle

   \section{Introduction}

   Kostka coefficients are important numbers appearing in many
   branches of mathematics, including representation theory, the
   theory of symmetric functions, and algebraic geometry (see,
   \emph{e.g.},~\cite{Ful97, Sta77, Sta99} and references
   therein).  Given a dominant weight \( \lambda \) and a weight
   \( \beta \) of the Lie algebra \( \gl_{r}(\C) \), the Kostka
   coefficient \( K_{\lambda \beta} \) is the dimension of the
   weight subspace with weight \( \beta \) of the irreducible
   representation \( V_{\lambda} \) of \( \gl_{r}(\C)
   \)~\cite{FH91}.  In the theory of symmetric functions, Kostka
   coefficients are defined by the expansion of Schur functions \(
   s_{\lambda} \) into monomials.  That is, given a partition \(
   \lambda \) of \( N \in \N \) with \( r \) parts, we have that
   \begin{equation}
   \label{symmfunc}
      s_{\lambda}
      =
      \sum_{%
         \substack{%
            \text{compositions} \\
            \text{\( \beta \) of \( N \)}
         }
      }
      K_{\lambda \beta} x^{\beta},
   \end{equation}
   where \( x^{\beta} = x_{1}^{\beta_{1}} x_{2}^{\beta_{2}} \dotsb
   x_{r}^{\beta_{r}} \).
   
   Since the parameters defining a Kostka coefficient are
   themselves vectors, they may be ``stretched'' by an integer
   scaling factor \( n \).  This procedure defines a function \( n
   \mapsto K_{n\lambda, n\beta} \), which, following~\cite{KTT04},
   we call a \emph{stretched Kostka coefficient}.  We denote this
   function by \( \K_{\lambda \beta}(n) = K_{n\lambda, n\beta} \).
   Kirillov and Reshetikhin have shown that \( \K_{\lambda
   \beta}(n) \) is a polynomial function of \( n \)~\cite{KR86}
   (see also \cite[Proposition 2.6]{DLM04}).  King, Tollu, and
   Toumazet have conjectured that these polynomials have only
   positive coefficients~\cite{KTT04}, and they have given a
   conjectural expression for the degree of \( \K_{\lambda
   \beta}(n) \)~\cite{KTT06}.  The main result of this note
   (Theorem \ref{main} below) is that the stretched Kostka
   coefficients indeed have the degrees conjectured
   in~\cite{KTT06}.

   For our purposes, it will suffice to assume that \( \lambda \)
   is a partition and \( \beta \) is a composition of the same
   length and size as \( \lambda \).  That is, we take \( \lambda
   = (\lambda_{1}, \dotsc, \lambda_{r}) \in \smash{\N^{r}} \) such
   that \( \lambda_{1} \geq \dotsb \geq \lambda_{r} \) and \(
   \beta = (\beta_{1}, \dotsc, \beta_{r}) \in \N^{r} \) such that
   \( \sum_{i} \beta_{i} = \sum_{i} \lambda_{i} \).  The Kostka
   coefficients \( K_{\lambda \beta} \) are indexed by such pairs
   \( (\lambda, \beta) \).  We sometimes write \( \lambda =
   (\kappa_{1}^{v_{1}}, \dotsc, \kappa_{m}^{v_{m}}) \) to indicate
   that \( \lambda \) has \( v_{p} \) parts equal to \( \kappa_{p}
   \) for \( 1 \leq p \leq m \).  The use of this notation always
   presumes that \( v_{p} \geq 1 \) for \( 1 \leq p \leq m \).
   
   Let \( \size{\lambda} = \sum_{i} \lambda_{i} \) denote the size
   of \( \lambda \).  Given an arbitrary sequence of nonnegative
   integers \( \beta = (\beta_{1}, \dotsc, \beta_{r}) \in \N^{r}
   \), let \( \prt{\beta} = (\prt{\beta}_{1}, \dotsc,
   \prt{\beta}_{r}) \) be the unique partition that may be
   produced by permuting the terms of \( \beta \).  We say that \(
   \lambda \) \emph{dominates} \( \beta \), denoted \( \beta
   \domby \lambda \), if \( \size{\lambda} = \size{\beta} \) and
   \( \sum_{k=1}^{i} \lambda_{k} \geq \sum_{k=1}^{i}
   \prt{\beta}_{k} \) for \( 1 \leq i < r \).  If, in addition, \(
   \sum_{k=1}^{i} \lambda_{k} > \sum_{k=1}^{i} \prt{\beta}_{k} \)
   for \( 1 \leq i < r \), we write \( \beta \primby \lambda \),
   and we say that \( \lambda \) and \( \beta \) form a
   \emph{primitive pair}.

   \begin{thm}[Proved on p.~\pageref{proof:GTdimproof}]
   \label{main}
      Suppose that \( \lambda = (\kappa_{1}^{v_{1}}, \dotsc,
      \kappa_{m}^{v_{m}}) \in \N^{r} \) is a partition and \(
      \beta \) is a composition such that \( \beta \primby \lambda
      \).  Then the degree of the stretched Kostka coefficient \(
      \K_{\lambda \beta}(n) \) is given by
      \begin{equation}
      \label{dimformula}
         \deg \K_{\lambda \beta}(n)
         =  \binom{r-1}{2}
            -
            \sum_{p=1}^{m} \binom{v_{p}}{2}
      \end{equation}
      \textup{(}where we evaluate \( \binom{1}{2} = 0
      \)\textup{)}%
      \footnote{%
         Equation \eqref{dimformula} was asserted in~\cite[Section
         7.5]{Kir01}.  However, the statement given there omits
         the condition that \( \lambda \) and \( \beta \) form a
         primitive pair.  Without this condition, the right-hand
         side of \eqref{dimformula} only gives an upper bound, as
         was observed in \cite{BGR04}.  For example, take \(
         \lambda = (4, 2, 1) \) and \( \beta = (3, 3, 1) \).  Then
         \( \deg \K_{\lambda \beta}(n) = 0 \), but the right-hand
         side of \eqref{dimformula}
         evaluates to \( 1 \).%
      }%
      .
   \end{thm}
   
   \begin{exm}
      Let \( \lambda = (4, 2, 2, 0, 0, 0) = (4^1, 2^2, 0^3) \) and
      \( \beta = (3, 1, 1, 1, 1, 1) \).  Then \( \lambda \primby
      \beta \), so Theorem \ref{main} applies.  In this case, we
      have \( r = 6 \), \( v_1 = 1 \), \( v_2 = 2 \), and \( v_3 =
      3 \), so the degree of \( \K_{\lambda \beta}(n) \) is
      \begin{equation*}
         \binom{6-1}{2}
         -  \binom{1}{2}
         -  \binom{2}{2}
         -  \binom{3}{2}
         =  6.
      \end{equation*}
   \end{exm}
   
   As stated, Theorem \ref{main} gives the degree of a stretched
   Kostka coefficient only when \( \lambda \) and \( \beta \) are
   a primitive pair.  However, Berenstein and Zelevinsky have
   shown that all Kostka coefficients factor into a product of
   Kostka coefficients indexed by primitive pairs~\cite{BZ90}.  It
   follows from this factorization that Theorem \ref{main}
   suffices to describe the degrees of stretched Kostka
   coefficients in all cases.

   The factorization of Kostka coefficients works as follows.  It
   is well known that \( K_{\lambda \beta} \) is invariant under
   permutations of the coordinates of \( \beta \).  For example,
   this follows from equation~\eqref{symmfunc} and the fact that
   Schur functions are symmetric.  Consequently, \( \K_{\lambda
   \beta}(n) = \K_{\lambda \prt{\beta}}(n) \).  In particular, to
   compute the degree of stretched Kostka coefficients, we need
   only consider the case where \( \beta \) is a partition.

   Suppose that \( \lambda, \beta \in \Z^{r} \) are both
   partitions with \( \size{\lambda} = \size{\beta} \).  If \(
   \lambda \) and \( \beta \) do not form a primitive pair, then
   we may write \( \lambda \) and \( \beta \), respectively, as
   concatenations of partitions such that each of the partitions
   contained in \( \lambda \) forms a primitive pair with the
   corresponding partition contained in \( \beta \).  More
   precisely, there exists a unique sequence of integers
   \begin{equation*}
      1 = i_{1} < i_{2} < \dotsb < i_{s} < i_{s+1} = r+1
   \end{equation*}
   such that each pair
   \begin{equation*} 
      \lambda^{(t)}
      =
      (\lambda_{i_{t}}, \dotsc, \lambda_{i_{t+1}-1})
      ,\,
      \beta^{(t)}
      =
      (\beta_{i_{t}}, \dotsc, \beta_{i_{t+1}-1}),
      \quad
      1 \le t \le s
   \end{equation*}
   is primitive.  We then have
   \begin{equation*} 
      K_{\lambda \beta} 
      =
      \prod_{t=1}^{s}
      K_{\lambda^{(t)} \beta^{(t)}}.
   \end{equation*}
   This observation of Berenstein and Zelevinsky~\cite{BZ90} is
   the justification for the terminology \emph{primitive pair}.
   Since the set of indices \( i_{1}, \dotsc, i_{s+1} \)
   decomposing \( (\lambda, \beta) \) into primitive pairs does
   not change when we scale \( \lambda \) and \( \beta \) by a
   parameter \( n \), this decomposition carries over to the
   stretched Kostka coefficients:
   \begin{equation*} 
      \K_{\lambda \beta}(n)
      =
      \prod_{t=1}^{s}
      \K_{\lambda^{(t)} \beta^{(t)}}(n).
   \end{equation*}
   Consequently, we have that
   \begin{equation*}
      \deg \K_{\lambda \beta}(n)
      =
      \sum_{t=1}^{s}
      \deg \K_{\lambda^{(t)} \beta^{(t)}}(n),
   \end{equation*}
   where each term in the sum may be computed using Theorem
   \ref{main}.
   
   This note is organized as follows.  In Section
   \ref{sec:dimensions}, we discuss the polyhedral model of Kostka
   coefficients using Gelfand--Tsetlin polytopes.  In this
   context, the degrees of stretched Kostka coefficients are the
   dimensions of these polytopes.  We prove that the expression in
   Theorem \ref{main} gives the dimension of a Gelfand--Tsetlin
   polytope under certain generic circumstances.  In Section
   \ref{sec:degKostka}, we show that these circumstances obtain
   precisely when the corresponding Kostka coefficient is indexed
   by a primitive pair, proving Theorem \ref{main}.

   \section{Dimensions of Gelfand--Tsetlin polytopes}
   \label{sec:dimensions}

   Several combinatorial interpretations of Kostka coefficients
   have appeared in the literature.  Most classically, \(
   K_{\lambda \beta} \) is the number of semi-standard Young
   tableaux with shape \( \lambda \) and content \( \beta \) (see,
   \emph{e.g.},~\cite{Sta97}).  Of particular interest for our
   study is the representation of Kostka coefficients as the
   number of lattice points in particular families of rational
   polytopes.  Gelfand and Tsetlin provided the first such
   model~\cite{GT50}, which we describe and employ in our study
   below.  
   
   The theory of lattice point enumeration has proved to be a
   powerful tool for understanding the behavior of Kostka
   coefficients.  For example, Billey, Guillemin, and Rassart have
   used vector partition functions to show that \( \K_{\lambda
   \beta}(n) \) can be expressed as a \emph{multivariate}
   piecewise polynomial in \( n \) and the cooridinates of \(
   \lambda \) and \( \beta \).  They also examined factorizations
   of these polynomials and gave upper bounds on their
   degrees~\cite{BGR04}.  More recently, King, Tollu, and Toumazet
   have introduced K-hive polytopes~\cite{KTT04}, which they used
   to motivate the conjectures mentioned in the introduction.
   Moreover, they deduce from their model additional information
   about the structure of the polynomials \( \K_{\lambda \beta}(n)
   \).  Among other results, they provide an interpretation for
   the roots of \( \K_{\lambda \beta}(n) \).

   Expressing Kostka coefficients using Gelfand--Tsetlin polytopes
   provides a natural geometric interpretation of the polynomial
   \( \K_{\lambda \beta}(n) \) and its degree.  To each Kostka
   coefficient \( K_{\lambda \beta} \), there corresponds a
   Gelfand--Tsetlin polytope \( GT_{\lambda \beta} \subset \R^{D}
   \) such that \( K_{\lambda \beta} = \card{GT_{\lambda \beta}
   \cap \Z^{D}} \), where \( D = \binom{r+1}{2} \).  As a
   consequence of the definition of Gelfand--Tsetlin polytopes
   (Definition \ref{GTpolydef} below), scaling \( \lambda \) and
   \( \beta \) by a positive integer \( n \) corresponds to
   dilating the polytope \( GT_{\lambda \beta} \) by \( n \):
   \begin{equation*} 
      \K_{\lambda \beta}(n) 
      =
      \card*{n GT_{\lambda \beta} \cap \Z^{D}}.
   \end{equation*}
   As mentioned in the introduction, Kirillov and Reshetikhin have
   shown that \( \K_{\lambda \beta}(n) \) is in fact a polynomial
   function of \( n \).  As is well known from the theory of
   lattice point enumeration in polyhedra, if the number of
   lattice points in an integer dilation \( n P \) of a polytope
   \( P \) is a polynomial function of \( n \), then the degree of
   that polynomial is the dimension of \( P \) (see,
   \emph{e.g.},~\cite[Theorem 4.6.25]{Sta97}).  Our proof of
   Theorem~\ref{main} depends on this interpretation of \( \deg
   \K_{\lambda \beta}(n) \) as the dimension of the
   Gelfand--Tsetlin polytope \( GT_{\lambda \beta} \), which we
   now define.

   Let \( X_{r} \) be the set of triangular arrays \(
   \smash{
   (x_{ij})_{1 \leq i \leq j \leq r}
   }
   \).  Note that \( X_r
   \) inherits a normed vector space structure under the obvious
   isomorphism \( X_r \cong 
   \smash{
   \R^{r(r+1)/2}
   }
   \).  Therefore, we will
   be able to speak of cones, lattices, and polytopes in \( X_{r}
   \).  It is customary to depict an array in \( X_{r} \) by
   arranging its entries as follows:
   \begin{equation*}
      \begin{matrix}
         x_{1r} &        & x_{2r} &        & x_{3r} & \multicolumn{3}{c}{\cdots} 
                                                                               & x_{rr} \\
                &        &        &        &        &        &        &        &        \\
                & \hdotsfor{7}                                                          \\
                &        &        &        &        &        &        &        &        \\
                &        & x_{13} &        & x_{23} &        & x_{33} &        &        \\
                &        &        &        &        &        &        &        &        \\
                & \phantom{x_{02}} 
                         &        & x_{12} &        & x_{22} &        & \phantom{x_{32}} 
                                                                               &        \\
                &        &        &        &        &        &        &        &        \\
                &        &        &        & x_{11} &        &        &        &
      \end{matrix}
   \end{equation*}
   Within the space \( X_{r} \), we define the cone of
   \emph{Gelfand--Tsetlin patterns}, or \emph{GT-patterns}, as
   follows.  A GT-pattern is a triangular array \( \gtp{x} =
   (x_{ij})_{1 \leq i \leq j \leq r} \in X_{r} \) satisfying the
   \emph{Gelfand--Tsetlin inequalities}
   \begin{equation}
   \label{eq:gtineq}
      x_{i,j+1} \ge x_{ij} \ge x_{i+1,j+1},
      \quad
      \text{for \( 1 \le i \le j\le r-1 \).}
   \end{equation}

   For any point \( \gtp{x} \in X_{r} \), we define the {\em
   highest weight} \( \hwt(\gtp{x}) \) of \( \gtp{x} \) to be the
   top row \( (x_{1r}, \dotsc, x_{rr}) \) of \( \gtp{x} \), and we
   define the {\em weight} \( \wt(\gtp{x}) = (\beta_{1}, \dotsc,
   \beta_{r}) \) of \( \gtp{x} \) by \( \beta_{1} = x_{11} \) and
   \( \beta_{j} = \sum_{i=1}^j x_{ij} - \sum_{i=1}^{j-1} x_{i,j-1}
   \) for \( 2 \leq j \leq r \).  Thus, \( \wt \) and \( \hwt \)
   are both linear maps \( X_{r} \rightarrow \R^{r} \).

   For each \( \lambda \in \R^{r} \), let \( GT_{\lambda} \) be
   the polytope of GT-patterns with highest weight \( \lambda \):
   \begin{equation*} 
      GT_{\lambda}
      =
      \setof{
         \gtp{x} \in X_{r}
         :
         \text{\( \gtp{x} \) is a GT-pattern and 
               \( \hwt(\gtp{x}) = \lambda \)}
      }.
   \end{equation*}
   Note that the GT inequalities~\eqref{eq:gtineq} force the top
   row of a GT-pattern to be weakly decreasing, so \( GT_{\lambda}
   = \emptyset \) if \( \lambda \) is not weakly decreasing.  (The
   converse is also true, since if \( \lambda \) is weakly
   decreasing, then \( (x_{ij})_{1 \leq i \leq j \leq r} \) with
   \( x_{ij} = \lambda_{i} \) is a GT-pattern in \( GT_{\lambda}
   \).)  For each \( \beta \in \R^{r} \), let \( W_{\beta} \subset
   X_{r} \) be the affine subspace of points in \( X_{r} \) with
   weight \( \beta \):
   \begin{equation*} 
      W_{\beta}
      =
      \setof{
         \gtp{x} \in X_{r}
         :
         \wt(\gtp{x}) = \beta
      }.
   \end{equation*}

   \begin{defn}
   \label{GTpolydef}
      Given \( \lambda, \beta \in \R^{r} \), the
      \emph{Gelfand--Tsetlin polytope} \( GT_{\lambda \beta} \) is
      the convex polytope of GT-patterns with highest weight \(
      \lambda \) and weight \( \beta \):
      \begin{equation*}
         GT_{\lambda \beta}
         =
         GT_{\lambda} \cap W_{\beta}.
      \end{equation*}
      Observe that each GT-pattern \( \gtp{x} \) lives in a unique
      GT-polytope.  We denote this polytope by \( GT(\gtp{x}) =
      GT_{\hwt(\gtp{x}), \wt(\gtp{x})}.  \)
   \end{defn}
    
   We associate to each GT-pattern a certain combinatorial object,
   called a \emph{tiling}.  Tilings of the GT-patterns in a
   GT-polytope have a poset structure isomorphic to the face
   lattice of the GT-polytope.  In particular, the tiling of a
   point \( \gtp{x} \) provides a straight-forward way to compute
   the dimension of the lowest-dimensional face of \( GT(\gtp{x})
   \) that contains \( \gtp{x} \) (Theorem \ref{facedimension}).

\begin{figure}
      \begin{center}
         \includegraphics{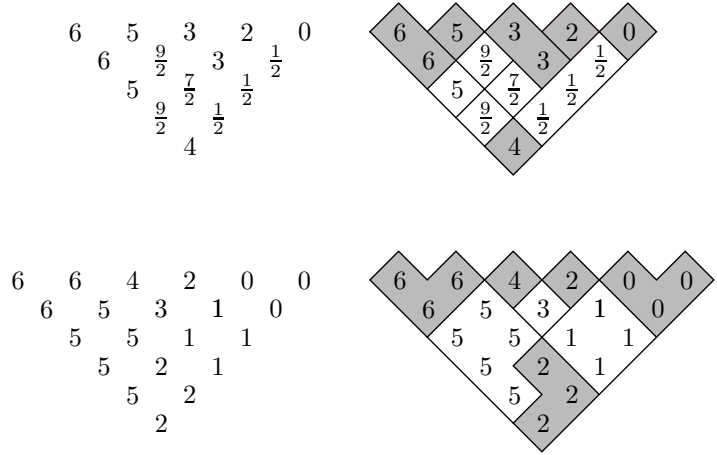}
         \caption{Tilings of GT-patterns.}
      \end{center}
      \label{fig:extiling}
   \end{figure}

   \begin{defn}
      The \emph{tiling} $\PP$ of $\gtp{x}$ is that partition of
      the set
      \begin{equation*}
         \mathcal{I}
         =
         \setof{(i,j) \in \Z^2 : 1 \leq i \leq j \leq r}
      \end{equation*}
      into subsets, called \emph{tiles}, that results from
      grouping together those entries in $\gtp{x}$ that are equal
      and adjacent.  More precisely, $\PP$ is that partition of \(
      \mathcal{I} \) such that two pairs $(i,j), (i',j')$ are in
      the same tile if and only if there are sequences
      \begin{gather*}
         i = i_1, i_2, \ldots, i_s = i' \\
         j = j_1, j_2, \ldots, j_s = j'
      \end{gather*}
      such that, for each $k \in \{1, \ldots, s-1 \}$, we have
      that
      \begin{equation*}
         (i_{k+1}, j_{k+1})
         \in
         \bigl\{
            (i_k + 1, j_k + 1),
            (i_k, j_k + 1),
            (i_k - 1, j_k - 1),
            (i_k, j_k - 1)
         \bigr\}
      \end{equation*}
      and $x_{i_{k+1} j_{k+1}} = x_{i_k j_k}$.
   \end{defn}
   In other words, the tiles are the connected components in the
   graph of a Gelfand--Tsetlin pattern.  See Figure
   \ref{fig:extiling} for examples of GT-patterns and their
   tilings.  The shading of some of the tiles in that figure is
   explained below.

   Given a GT-pattern \( \gtp{x} \) with tiling \( \PP \), we
   associate to \( \PP \) (or, equivalently, to \( \gtp{x} \)) a
   matrix \( \A_{\PP} \) (or \( A_{\gtp{x}} \)) as follows.
   Define the \emph{free rows} of \( \gtp{x} \) to be those that
   are neither the top nor the bottom row of \( \gtp{x} \).  The
   \emph{free tiles} \( T_{1}, \dotsc, T_{s} \) of \( \PP \) are
   those tiles in \( \PP \) that intersect only free rows of \(
   \gtp{x} \)---\textit{i.e.}, those tiles that contain neither \(
   (1,1) \) nor \( (i,r) \) for \( 1 \leq i \leq r \).  The
   remaining tiles are the \emph{non-free tiles}.  The order in
   which the free tiles are indexed will not matter for our
   purposes, but, for concreteness, we adopt the convention of
   indexing the free tiles in the order that they are initially
   encountered as the entries of \( \gtp{x} \) are read from left
   to right and bottom to top.  Define the \emph{tiling matrix} \(
   A_{\PP} = A_{\gtp{x}} = (a_{jk})_{1 \leq j \leq r-2, \, 1 \leq
   k \leq s} \) by
   \begin{equation*}
      a_{jk}%
      =%
      \# \setof{i : (i,j+1) \in T_k }.
   \end{equation*}
   That is, \( a_{jk} \) counts the number of entries in the \( j
   \)th free row of \( \gtp{x} \) that are contained in the free
   tile \( P_k \).  While a different choice of order for the free
   tiles would result in a tiling matrix with permuted columns,
   this will be immaterial because we will ultimately be
   interested only in the rank of \( \A_{\PP} \).  The terminology
   \emph{free tile} represents the fact that the entries in the
   non-free tiles of \( \gtp{x} \) are fixed once the weight and
   highest weight of \( \gtp{x} \) are fixed.

   \begin{exm}
      Two GT-patterns and their tilings are given in Figure
      \ref{fig:extiling}.  The unshaded tiles are the free tiles.
      Index the free tiles in the order that they are initially
      encountered as the entries of $\gtp{x}$ are read from left
      to right and bottom to top.  Then the associated tiling
      matrices are respectively
      \begin{equation*}
         \begin{bmatrix}
            1 & 1 & 0 & 0 & 0 \\
            0 & 1 & 1 & 1 & 0 \\
            0 & 1 & 0 & 0 & 1
         \end{bmatrix}
         \quad \text{and} \quad
         \begin{bmatrix}
            1 & 0 & 0 \\
            1 & 1 & 0 \\
            2 & 2 & 0 \\
            1 & 1 & 1 
         \end{bmatrix}.
      \end{equation*}
   \end{exm}

   % \pdfsyncstop
   % \pdfsyncstart

   The primary motivation for introducing tilings of GT-patterns
   is the following Theorem:
   
   \begin{thm}[\protect{\cite[Theorem 1.5]{DLM04}}]
   \label{facedimension}
      Suppose that $\PP$ is the tiling of a GT-pattern $\gtp{x}$.
      Then the dimension of the kernel of $A_\PP$ is equal to the
      dimension of the minimal \textup{(}dimensional\textup{)}
      face of the GT-polytope containing $\gtp{x}$.
   \end{thm}

   Theorem \ref{facedimension} was used in~\cite{DLM04} to study
   the properties of vertices of GT-polytopes, establishing in
   particular that they can have arbitrarily large denominators.
   In the present note, we move to the opposite end of the face
   lattice and apply the tiling machinery to points in the
   interior of GT-polytopes.  We use the notation \( \intr P \)
   to denote the \emph{relative interior} of a polytope \( P
   \)---that is, the interior of \( P \) with respect to the
   affine space that it spans.  As an immediate corollary to
   Theorem \ref{facedimension}, we get:
   \begin{cor}
   \label{intpointdim}
      If \( \gtp{x} \in \intr(GT_{\lambda \beta}) \), then \( \dim
      GT_{\lambda \beta} = \dim \ker A_{\gtp{x}} \).
   \end{cor}

\begin{figure}
      \center
      \includegraphics{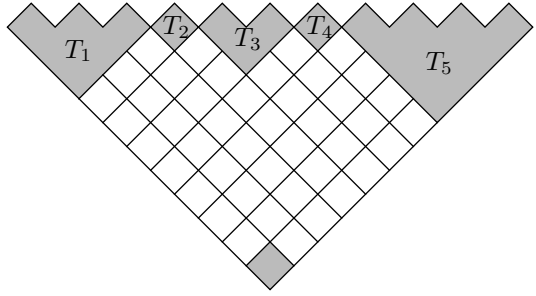}
         \caption{%
            The tiling of a point in the interior of \(
            GT_{\lambda} \).
         }%
      \label{fig:tiling}
   \end{figure}

   If \( GT_{\lambda} \ne \emptyset \), the tilings of GT-patterns
   in the relative interior of \( GT_{\lambda} \) have an easy
   characterization.  Observe that if a block of entries \(
   x_{kr},\, x_{k+1,r},\, \dotsc,\, x_{\ell r} \) in the top row
   of a GT-pattern all have the same value \( \kappa_{p} \), then
   the GT inequalities~\eqref{eq:gtineq} require that the entries
   \( x_{ij} \) with \( k \leq i \leq \ell-(r-j) \) all assume
   that same value \( \kappa_{p} \).  These are the entries that
   lie within the triangular region whose horizontal edge consists
   of the entries \( x_{kr},\, x_{k+1,r},\, \dotsc,\, x_{\ell r}
   \) and whose diagonal edges run parallel to the diagonals of
   the GT-pattern.  Let \( T_{p} \) be the tile of entries in this
   region.  That is, if \( \lambda = (\kappa_{1}^{v_{1}}, \dotsc,
   \kappa_{m}^{v_{m}}) \in \R^{r} \), put
   \begin{equation*} 
      T_{p}
      =
      \setof{
         (i,j) \in \mathcal{I}
         \,:\,
         \sum_{q=1}^{p-1} v_{q}
         <
         i
         \leq
         \sum_{q=1}^{p} v_{q} - (r - j)
      },
      \quad
      1 \leq p \leq m.
   \end{equation*}
   Define the \emph{generic interior tiling} \( \PP_{\lambda} \)
   associated with \( \lambda \) to be the tiling consisting of
   these \( T_{p} \)'s together with a distinct tile for each
   entry not contained in one of the \( T_{p} \)'s.  See Figure
   \ref{fig:tiling} for an example of a generic interior tiling
   when \( (v_{1}, \dotsc, v_{5}) = (3,1,2,1,4) \).  The shaded
   tile at the bottom of the pattern, and each of the unshaded
   free tiles, contains only a single entry.

   % \pdfsyncstop
   % \pdfsyncstart

   The next theorem characterizes when a GT-pattern \( \gtp{x} \)
   has the generic interior tiling and gives the dimension of the
   GT-polytope containing \( \gtp{x} \).  This is the main result
   from which Theorem \ref{main} will follow.

\begin{figure}
      \center
      \includegraphics{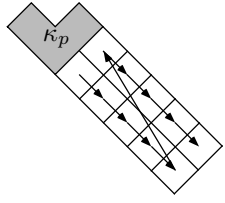}
      \caption{%
         Filling of a portion of the generic interior tiling.%
      }%
      \label{fig:filling}
   \end{figure}

   \begin{thm}
   \label{maincor}
      Let \( \lambda = (\kappa_{1}^{v_{1}}, \dotsc,
      \kappa_{m}^{v_{m}}) \in \R^{r} \).
      \begin{enumerate}
         \item
         A GT-pattern \( \gtp{x} \in GT_{\lambda} \) is in the
         relative interior of \( GT_{\lambda} \) if and only if
         the tiling of \( \gtp{x} \) is the generic interior
         tiling \( \PP_{\lambda} \).
         
         \item
         If \( \gtp{x} \) has the generic interior tiling \(
         \PP_{\lambda} \) and \( m \geq 2 \), then
         \begin{equation*} 
            \dim GT(\gtp{x})
            =
            \binom{r-1}{2}
            -
            \sum_{p=1}^{m} \binom{v_{p}}{2},
         \end{equation*}
         \textup{(}where we evaluate \( \binom{1}{2} = 0
         \)\textup{)}.
      \end{enumerate}
   \end{thm}
   \begin{proof}
      To prove part (1), note that the GT
      inequalities~\eqref{eq:gtineq} include the facet-defining
      inequalities of \( GT_{\lambda} \) in \( X_{r} \).  A point
      \( \gtp{x} \) lies in the relative interior of \(
      GT_{\lambda} \) if and only if \( \gtp{x} \) satisfies with
      equality \emph{only} those GT-inequalities that are
      satisfied with equality by every point in \( GT_{\lambda}
      \).
      
      We claim that the GT-inequalities satisfied with equality by
      every point in \( GT_{\lambda} \) are precisely the ones
      implied by forcing each entry in \( T_{p} \) to be equal to
      \( \kappa_{p} \) for \( 1 \leq p \leq m \).  Since every
      point in \( GT_{\lambda} \) must satisfy those equalities,
      it remains only to exhibit a point \( \gtp{x} \in
      GT_{\lambda} \) satisfying \emph{only} those equalities.
      That is, we need to construct an \( \gtp{x} \in GT_{\lambda}
      \) with the generic interior tiling.
      
      For \( 1 \leq p \leq m-1 \), consider the
      upper-left-to-lower-right diagonals that abut the tile \(
      T_{p} \).  Let \( S_{p} \) be the set of entries in these
      diagonals.  Apply a total order to \( S_{p} \) by reading
      each diagonal from left to right, and then reading the
      diagonals themselves from left to right.  Finally, construct
      \( \gtp{x} \) by assigning values to the entries in \( S_{p}
      \) according to a strictly monotonically decreasing function
      mapping \( S_{p} \) into the open interval \(
      \opint{\kappa_{p+1}, \kappa_{p}} \).  For example, in the
      portion depicted in Figure \ref{fig:filling}, we fill the
      tiles with strictly decreasing values in the order indicated
      by the arrows.  It is easy to see that \( \gtp{x} \) is a
      GT-pattern and that its tiling is \( \PP_{\lambda} \), so
      part (1) is proved.
      
      To prove part (2), suppose that \( \gtp{x} \) has the
      generic interior tiling.  Then, by part (1), we have that \(
      \gtp{x} \in \intr(GT(\gtp{x})) \).  Thus, we can apply
      Corollary \ref{intpointdim} to compute \( \dim GT(\gtp{x})
      \).  The hypothesis that \( m \geq 2 \) implies that every
      free row of \( \gtp{x} \) contains a free tile.  Hence, each
      of the \( r-2 \) rows of the tiling matrix \( A_{\gtp{x}} \)
      contains a nonzero entry.  Moreover, since every free tile
      of the generic interior tiling \( \PP_{\lambda} \) contains
      only a single entry, every column of \( \A_{\gtp{x}} \)
      contains only a single 1.  Thus, \( \A_{\gtp{x}} \) is in
      reduced row echelon form (perhaps after a suitable
      permutation of its columns, which amounts to re-indexing the
      free tiles).  This means that
      \begin{align*} 
         \dim GT(\gtp{x})
         &= \dim \ker \A_{\gtp{x}} \\
         &= (\text{%
               \# of columns of \( \A_{\gtp{x}} \)%
            })
            -
            (\text{%
               dimension of row span of \( \A_{\gtp{x}} \)%
            }) \\
         &= (\text{%
               \# of free tiles in \( \PP_{\lambda} \)%
            })
            -                
            (r - 2).
      \end{align*}
      Thus the computation reduces to finding the number of free
      tiles in \( \PP_{\lambda} \), which is easily done:
      \begin{equation*} 
         \binom{r}{2} 
         -
         \sum_{p=1}^{m} \bigl(\card*{T_{p}} - v_{p}\bigr) - 1 
         =
         \binom{r}{2}
         -
         \sum_{p=1}^{m} \binom{v_{p}}{2} - 1.
      \end{equation*}
      Hence,
      \begin{equation*} 
         \dim GT(\gtp{x})
         =
         \binom{r}{2}
         -
         \sum_{p=1}^{m} \binom{v_{p}}{2} - (r - 1)
         =
         \binom{r-1}{2}
         -
         \sum_{p=1}^{m} \binom{v_{p}}{2},
      \end{equation*}
      as claimed.
   \end{proof}

   % \pdfsyncstop
   % \pdfsyncstart

   \section{The degree of stretched Kostka coefficients}
   \label{sec:degKostka}

   In the previous section, we showed that if the interior points
   of \( GT_{\lambda \beta} \) have the generic interior tiling,
   then the dimension of \( GT_{\lambda \beta} \) is given by the
   expression in Theorem \ref{main}.  To complete our proof of
   Theorem \ref{main}, it remains only to show that, if \( \beta
   \primby \lambda \), then the interior points of \( GT_{\lambda
   \beta} \) have the generic interior tiling.  To this end, we
   call upon a well-known fact from the theory of convex
   polytopes.

   \begin{lem}
   \label{inttoint}
      Given a polytope \( P \subset \R^{m} \) and a linear map \(
      \pi \colon \R^{m} \to \R^{n} \), we have that \(
      \intr(\pi(P)) \subseteq \pi(\intr{P}) \).
   \end{lem}

   The converse containment is also true and easy to prove, but it
   is unnecessary for our purposes.  It is also worth mentioning
   that this is the point at which the geometry of \emph{convex}
   polyhedra is crucial to the argument.  Lemma \ref{inttoint}
   does not hold for every polyhedral set, even if it is connected
   and full-dimensional.  For completeness, we now give a proof of
   Lemma \ref{inttoint}.

   \begin{proof}[Proof of Lemma \ref{inttoint}]
      Without loss of generality, we assume that \( \dim P = m \)
      and \( \dim \pi(P) = n \).  Suppose that \( y \notin
      \pi(\intr{P}) \).  We show that \( y \notin \intr(\pi(P)) \)
      by exhibiting an ``exit vector'' \( \bar{y} \) such that,
      for every \( \epsilon > 0 \), \( y + \epsilon \bar{y} \notin
      \pi(P) \).  Since we assume that \( \pi(P) \) is
      full-dimensional, this will prove the claim.
      
      Since \( y \notin \pi(\intr P) \), we have that \(
      \pi^{-1}(y) \cap \intr P = \emptyset \).  The sets \(
      \pi^{-1}(y) \) and \( \intr P \) are both convex, so there
      is a hyperplane \( H \) separating them.  Let \( \bar{x} \in
      \R^{m} \) be a normal to \( H \) pointing away from \( P \).
      Then, for every \( x' \in \pi^{-1}(y) \) and \( \epsilon > 0
      \), we have that \( x' + \epsilon \bar{x} \notin P \).  Let
      \( \bar{y} = \pi(\bar{x}) \).  Suppose that \( \epsilon > 0
      \) and that \( x \in \R^{m} \) is such that \( \pi(x) = y +
      \epsilon \bar{y} \).  Then \( x - \epsilon \bar{x} \in
      \pi^{-1}(y) \), so \( x = x - \epsilon \bar{x} + \epsilon
      \bar{x} \notin P \).  In other words, \( y + \epsilon
      \bar{y} \notin \pi(P) \), proving the claim.
   \end{proof}

   Let \( \Dom(\lambda) \subset \R^{r} \) be the image of \(
   GT_{\lambda} \) under the map \( \wt \colon X_{r} \to \R^{r}
   \).  Note that for \( \beta \in \Dom(\lambda) \), we have \(
   GT_{\lambda \beta} = \wt^{-1}(\beta) \cap GT_{\lambda} \).  It
   is well known that \( GT_{\lambda \beta} \ne \emptyset \) if
   and only if \( \beta \domby \lambda \)~\cite{KB95}.  For, if \(
   \beta \domby \lambda \), then \( K_{\lambda \beta} > 0 \) (see,
   \textit{e.g.},~\cite[Exercise A.11]{FH91}), so \( GT_{\lambda
   \beta} \) contains an integral point.  Conversely, if \(
   GT_{\lambda \beta} \ne \emptyset \), then some integral
   multiple \( n GT_{\lambda \beta} = GT_{n\lambda, n\beta} \)
   contains an integral point.  Hence, \( K_{n\lambda, n\beta} > 0
   \), so \( n\beta \domby n\lambda \).  Since the relative order
   of two weights is not changed by scaling both by \( n \), it
   follows that \( \beta \domby \lambda \).  This establishes the
   following lemma:

   \begin{lem}
   \label{permint}
      Suppose that \( \lambda \in \Z^{r} \) is a partition.  Then
      \( \Dom(\lambda) = \setof{ \beta \in \R^{r} : \beta \domby
      \lambda } \).  Consequently, if \( \beta \primby \lambda \),
      then \( \beta \in \intr(\Dom(\lambda)) \).
   \end{lem}

   Putting together the preceding results, we are now ready to
   prove Theorem~\ref{main} from page~\pageref{main}.

   \begin{proof}[Proof of Theorem~\ref{main}]
   \label{proof:GTdimproof}
      The claim is trivial if \( r=1 \), so suppose that \( r \ge
      2 \).  From Lemmas~\ref{permint} and~\ref{inttoint}, we have
      that
      \begin{equation*}
         \beta 
         \in
         \intr(\Dom(\lambda))
         \subset
         \wt(\intr(GT_{\lambda})).
      \end{equation*}
      Hence,
      \begin{equation*}
         GT_{\lambda \beta} \cap \intr(GT_{\lambda})
         =
         \wt^{-1}(\beta) \cap \intr(GT_{\lambda})
         \ne
         \emptyset.
      \end{equation*}
      Choose \( \gtp{x} \in GT_{\lambda \beta} \cap
      \intr(GT_{\lambda}) \).  By part (1) of
      Theorem~\ref{maincor}, \( \gtp{x} \) has the generic
      interior tiling \( \PP_{\lambda} \).  Note that since \(
      \beta \primby \lambda \), we must have that \( m \ge 2 \),
      for otherwise the GT-inequalities would force \( \beta_{1} =
      \lambda_{1} \).  Hence, we may apply part (2) of
      Theorem~\ref{maincor} to compute the dimension of \(
      GT(\gtp{x}) = GT_{\lambda \beta} \).  Since the dimension of
      \( GT_{\lambda \beta} \) is the degree of \( \K_{\lambda
      \beta}(n) \)~\cite[Theorem 4.6.25]{Sta97}, this yields
      \begin{align*}
         \deg \K_{\lambda \beta}(n)
         &= \dim GT_{\lambda \beta} = \dim GT(\gtp{x}) \\
         &= \binom{r-1}{2}
            -
            \sum_{p=1}^{m} \binom{v_{p}}{2},
      \end{align*}
      as claimed.
   \end{proof}

   \section*{Acknowledgements}

   The author wishes to thank Christophe Tollu and Frederic
   Toumazet for exposing him to their conjecture, and Jes\'us De
   Loera for helpful comments and suggestions.  The author also
   thanks the anonymous referees for suggestions improving the
   exposition.

%    \bibliographystyle{hamsplainf}
%    \bibliography{MasterBibliography}

\begin{thebibliography}{10}

\bibitem{BZ90}
A.~D. Berenstein and A.~V. Zelevinsky, \emph{When is the multiplicity of a
  weight equal to {$1$}?}, Funktsional. Anal. i Prilozhen. \textbf{24} (1990),
  no.~4, 1--13, 96.

\bibitem{BGR04}
S.~Billey, V.~Guillemin, and E.~Rassart, \emph{A vector partition function for
  the multiplicities of {$\germ{sl}\sb k\mathbb{C}$}}, J. Algebra \textbf{278}
  (2004), no.~1, 251--293.

\bibitem{DLM04}
J.~A. De~Loera and T.~B. McAllister, \emph{Vertices of {G}elfand-{T}setlin
  polytopes}, Discrete Comput. Geom. \textbf{32} (2004), no.~4, 459--470,
  \mbox{arXiv:math.CO/0309329}.

\bibitem{Ful97}
W.~Fulton, \emph{Young tableaux}, London Mathematical Society Student Texts,
  vol.~35, Cambridge University Press, Cambridge, 1997, With applications to
  representation theory and geometry.

\bibitem{FH91}
W.~Fulton and J.~Harris, \emph{Representation theory}, Graduate Texts in
  Mathematics, vol. 129, Springer-Verlag, New York, 1991, A first course,
  Readings in Mathematics.

\bibitem{GT50}
I.~M. Gelfand and M.~L. Tsetlin, \emph{Finite-dimensional representations of
  the group of unimodular matrices}, Doklady Akad. Nauk SSSR (N.S.) \textbf{71}
  (1950), 825--828, English translation in \emph{Collected Papers, Vol. II},
  Springer, Berlin, 1988, 653--656.

\bibitem{KTT04}
R.~C. King, C.~Tollu, and F.~Toumazet, \emph{Stretched
  {L}ittlewood-{R}ichardson and {K}ostka coefficients}, Symmetry in physics,
  CRM Proc. Lecture Notes, vol.~34, Amer. Math. Soc., Providence, RI, 2004,
  pp.~99--112.

\bibitem{KTT06}
\bysame, \emph{The hive model and the factorisation of {K}ostka coefficients},
  S{\'e}minaire Lotharingien de Combinatoire \textbf{54A} (2006), Article
  B54Ah.

\bibitem{Kir01}
A.~N. Kirillov, \emph{Ubiquity of {K}ostka polynomials}, Physics and
  combinatorics 1999 (Nagoya), World Sci. Publishing, River Edge, NJ, 2001,
  \mbox{arXiv:math.QA/9912094}, pp.~85--200.

\bibitem{KB95}
A.~N. Kirillov and A.~D. Berenstein, \emph{Groups generated by involutions,
  {G}el\cprime fand-{T}setlin patterns, and combinatorics of {Y}oung tableaux},
  Algebra i Analiz \textbf{7} (1995), no.~1, 92--152, {\em translation in St.
  Petersburg Math. J.} \textbf{7} (1996), no. 1, 77--127.

\bibitem{KR86}
A.~N. Kirillov and N.~Y. Reshetikhin, \emph{The {B}ethe ansatz and the
  combinatorics of {Y}oung tableaux}, Zap. Nauchn. Sem. Leningrad. Otdel. Mat.
  Inst. Steklov. (LOMI) \textbf{155} (1986), no.~Differentsialnaya Geometriya,
  Gruppy Li i Mekh. VIII, 65--115, 194, {\em translation in J. Soviet Math.}
  \textbf{41} (1988), no. 2, 925--955.

\bibitem{Sta77}
R.~P. Stanley, \emph{Some combinatorial aspects of the {S}chubert calculus},
  Combinatoire et repr\'esentation du groupe sym\'etrique (Actes Table Ronde
  CNRS, Univ. Louis-Pasteur Strasbourg, Strasbourg, 1976), Springer, Berlin,
  1977, pp.~217--251. Lecture Notes in Math., Vol. 579.

\bibitem{Sta97}
\bysame, \emph{Enumerative combinatorics. {V}ol. 1}, Cambridge Studies in
  Advanced Mathematics, vol.~49, Cambridge University Press, Cambridge, 1997,
  With a foreword by Gian-Carlo Rota, Corrected reprint of the 1986 original.

\bibitem{Sta99}
\bysame, \emph{Enumerative combinatorics. {V}ol. 2}, Cambridge Studies in
  Advanced Mathematics, vol.~62, Cambridge University Press, Cambridge, 1999,
  With a foreword by Gian-Carlo Rota and appendix 1 by Sergey Fomin.

\end{thebibliography}
   
   \def\cprime{$'$}
\providecommand{\bysame}{\leavevmode\hbox to3em{\hrulefill}\thinspace}
\providecommand{\href}[2]{#2}

\end{document}